\documentclass[12pt,leqno]{article}

\usepackage[cp1250]{inputenc}
\usepackage[OT4]{fontenc}
\usepackage{amsthm}
\usepackage{amsfonts}
\usepackage{amsmath}
\usepackage{amssymb}

\newcommand{\R}{\mathbb{R}}
\newcommand{\C}{\mathbb{C}}

\newcommand{\bo}{{\bf 0}}
\newcommand{\bO}{{\bf 0}}

\newcommand{\Opis}{{\cal O}}

\newtheorem{theorem} {Theorem}[section]
\newtheorem{theorem*}{Theorem}
\newtheorem{prop*} {Proposition}
\newtheorem{lemma*}{Lemma}
\newtheorem{lemma}[theorem]{Lemma}

\newtheorem{cor*}{Corollary}
\newtheorem{prop}[theorem] {Proposition}

\newtheorem{definition*}{Definition}
\newtheorem{rem}[theorem]{Remark}

\theoremstyle{definition}

\begin{document} 
 

\begin{center}
{\Large On bifurcation of cusps}\\[1em]
by  Zbigniew Szafraniec\\[1em]
\end{center}
{\bf Abstract} Let $f={\bf R}\times{\bf R}^2,\bo\rightarrow {\bf R}^2,\bo$ be an analytic mapping
having a critical point at the origin. There is the corresponding
one-parameter family of mappings $f_t=f(t,\cdot):{\bf R}^2\rightarrow{\bf R}^2$.

There will be presented effective algebraic methods of computing  the number
of cusps of $f_t$, where $0<|t|\ll 1$, emanating from the origin and having a positive/negative cusp degree..\\[1em]


\section{Introduction}\label{wstep}
Mappings between surfaces are a natural object of study
in the theory of singularities.  Whitney \cite{whitney} proved
that critical points of such a generic mapping are folds and cusps.
There are several results concerning relations between the topology of surfaces
and the topology of the critical locus of a mapping (see \cite{levine}, \cite{thom}, \cite{whitney}).
Singularities of map germs of the plane into the plane were studied in
\cite{fukudaishikawa}, \cite{fukuda2}, \cite{gaffneymond},  \cite{moyaballesteros1},  
\cite{moyaballesteros2}, \cite{rieger}.

Let $f_t$, where $t\in\R$, be an analytic family of plane-to-plane mappings
with $f_0$ having a critical point at the origin.
Under some natural assumptions  there is a finite family
of cusp points of $f_t$ bifurcating from the origin.
There are important results \cite[Theorem 3.1]{fukudaishikawa}, \cite[Proposition 7.1]{moyaballesteros2}
concerning the parity of the number of those points.

In this paper we show how to compute the number of cusps of $f_t$ which are
represented by germs having either positive or negative local topological degree
(see Theorem \ref{rodziny6}).

The paper is organized as follows. In Sections \ref{powierzchnie} and \ref{kielki},
we collect some useful facts.
The curve in $\R\times\R^2$ consisting of points $(t,x)$, where $x$ is a cusp point of $f_t$, is defined by three analytic equations,
so that  it is not a complete intersection.
In  Section \ref{liczbagalezi} we show how to adopt in this case some  more general techniques 
from \cite{nowelszafraniec}
concerning curves in $\R^n$ defined by $m$ equations, where $m\geq n$.

In Sections \ref{owale} and \ref{rodziny}, we prove the main result.
In Section \ref{przyklad} we present  examples computed by a computer. We have implemented our algorithm
with the help of {\sc Singular} \cite{singular}. We have also used a computer program written by
{\L}\c{e}cki \cite{leckiszafraniec1}.


\section{Mappings between surfaces}\label{powierzchnie}

Let $(M,\partial M)$ and $(N,\partial N)$ be  compact oriented connected surfaces, and let 
$f:M\rightarrow N$ be a smooth mapping such that $f^{-1}(\partial N)=\partial M$.
Assume that
\begin{itemize}
\item [(i)]  every point in $M$ is either a fold point, a cusp point or a regular point,
and there is only a finite number of cusps which all belong to $M\setminus\partial M$,
\item [(ii)] the $1$-dimensional manifold consisting of fold points is transverse to $\partial M$,
so that $f|\partial  M:\partial M\rightarrow \partial N$ is locally stable, i.e. its critical points
are non-degenerate.
\end{itemize}

We shall write $M^-$ for the closure in $M$ of the set of regular points at which
$f$ does reverse the orientation.

If  $p\in M\setminus\partial M$ is a cusp point, we define $\mu(p)$
to be the local topological degree of the germ $f:(M,p)\rightarrow (N,f(p))$.
Put
$$\operatorname{cusp\  deg}\,(f)=\sum \mu(p),$$
where $p$ runs through the set of all cusp points of $f$.

Fukuda and Ishikawa \cite{fukudaishikawa} have generalized the result by Quine \cite{quine}
concerning surfaces withouth boundary, proving

\begin{theorem}\label{powierzchnie1}
Let $M, N$ and $f$ be as above and $\partial M\neq\emptyset$. Then
$$\operatorname{cusp\ deg}\,( f)=2\chi(M^-)+(\deg\, f|\partial M)\chi(N)-\chi(M)-\#C(f|\partial M)/2,$$
where  $C(f|\partial M)$ is the set of critical points of $f|\partial M$.
\end{theorem}

In fakt, in \cite{fukudaishikawa} there is a stronger assumption that both
$f:M\rightarrow N$ and $f|\partial M:\partial M\rightarrow\partial N$ are $C^\infty$--stable mappings.
However, if $f$ satisfies (i), (ii) then there exists its $C^\infty$--stable perturbation $\tilde f$,
which is arbitrary close to $f$ in $C^\infty$--Whitney topology,
such that all corresponding numbers associated to $f$ and $\tilde f$ which appear in the above theorem
stay the same.

Let $f=(f_1,f_2):U\rightarrow\R^2$, where $U\subset\R^2$ is open, be a smooth mapping. Set
$J=\partial(f_1,f_2)/\partial(x_1,x_2),\ F_i=\partial(f_i,J)/\partial(x_1,x_2),\  i=1,2.$
Applying the same arguments as in the proof of  \cite[Proposition 2, p. 815]{krzyzanowskaszafraniec}
one gets

\begin{prop}\label{powierzchnie2}
The set of all  common solutions in $U$ of the system of equations
  $J=F_1=F_2=\partial(F_1,J)/\partial(x_1,x_2)=\partial(F_2,J)/\partial(x_1,x_2)=0$
is empty if and only if the set of critical points of $f$ consists of either fold or cusp points.

If that is the case then the set of cusp points is discrete and equals $\{J=F_1=F_2=0\}$.
\end{prop}

\section{Families of germs}\label{kielki}

In this section we  recall some useful facts concerning
1-parameter families of real analytic germs.

For $r>0$, let 
$D^n(r)=\{x\in\R^n\, |\, \parallel x\parallel \leq r\}$, and
$S^{n-1}(r)=\partial D^n(r)$.
We shall write $(t,x)=(t,x_1,\ldots,x_n)\in\R\times\R^n$.
Assume $J(t,x):\R\times\R^n,\bO\rightarrow \R,0$ is an analytic function defined in a neighbourhood of the origin having a critical point at $\bo$. We shall write
$$L_0=\{x\in S^{n-1}(r)\, |\, J(0,x)= 0\},$$
$$M_t^{-}=\{x\in D^n(r)\, |\, J( t,x)\leq 0\},$$
where $0<|t|\ll r\ll 1$.

Let $f:\R\times\R^n,\bo\rightarrow\R^n,\bo$ be an analytic mapping. Put $f_t(x)=f(t,x)$.
Suppose that there exists a small $r>0$ such that $f_0^{-1}(\bo)\cap D^n(r)=\{\bo\}$.
For $0<\delta\ll r$, put $\tilde{S}_t^{n-1}(\delta)=f_t^{-1}(S^{n-1}(\delta))\cap D^n(r)$
and $\tilde{D}_t^n(\delta)=f_t^{-1}(D^n(\delta))\cap D^n(r)$.
We shall write
$$\tilde{L}_0=\{x\in \tilde{S}_0^{n-1}(\delta)\, |\, J(0,x)= 0\},$$
$$\tilde{M}_t^{-}=\{x\in \tilde{D}_t^n(\delta)\, |\, J( t,x)\leq 0\},$$
where $0<|t|\ll \delta\ll 1$.

\begin{lemma}\label{powtorka2}
We have $\chi(\tilde{M}_t^{-})=\chi(M_t^{-})$ and $\chi(\tilde{L}_0  )=\chi(L_0)$.
\end{lemma}
\noindent{\em Proof.} There exist small positive $\delta_1<\delta_2$, $ r_1<r_2$ 
and $t_0$, such that for $0<|t|<t_0$ we have
$$\{ x\in\tilde{D}_t^n(\delta_1)\, |\,  J( t,x)\leq 0  \} \subset \{  x\in D(r_1)\, |\,  J( t,x)\leq 0  \}$$
$$\subset\{ x\in\tilde{D}_t^n(\delta_2)\, |\, J( t,x)\leq 0  \} \subset \{  x\in D(r_2)\, |\,  J( t,x)\leq 0  \},$$
and inclusions
$$\{ x\in\tilde{D}_t^n(\delta_1)\, |\,  J( t,x)\leq 0  \} \subset \{ x\in\tilde{D}_t^n(\delta_2)\, |\,  J( t,x)\leq 0  \},$$
$$\{  x\in D(r_1)\, |\, J( t,x)\leq 0  \} \subset \{  x\in D(r_2)\, |\,  J( t,x)\leq 0  \}$$
induce isomorphisms of corresponding homology groups. Then
$$\chi (\tilde{M}_t^{-})= \chi(\{ x\in\tilde{D}_t^n(\delta_1)\, |\,  J( t,x)\leq 0  \} )$$
$$=\chi(\{  x\in D(r_2)\, |\,  J( t,x)\leq 0  \})=\chi(M_t^{-}).$$
The proof of the second assertion is similar. $\Box$

Define a mapping $d_0:\R^n,\bo\rightarrow \R^n,\bo$ by
$$d_0(x)=\left(  \frac{\partial J}{\partial x_1}(0,x),\ldots,  \frac{\partial J}{\partial x_n}(0,x)  \right),$$
and mappings $d_1,d_2:\R\times\R^n,\bo\rightarrow\R\times\R^n,\bo$, by
$$d_1(t,x)=\left(    \frac{\partial J}{\partial t}(t,x), \frac{\partial J}{\partial x_1}(t,x),\ldots, \frac{\partial J}{\partial x_n}(t,x)       \right),$$
$$d_2(t,x)=\left(    J(t,x), \frac{\partial J}{\partial x_1}(t,x),\ldots, \frac{\partial J}{\partial x_n}(t,x)       \right).$$
Applying directly results by Fukui \cite{fukui2} and Khimshiasvili \cite{khimshiashvili1, khimshiashvili2} we get

\begin{theorem}\label{powtorka1}
Suppose that the origin is isolated in $d_0^{-1}(\bo)$, $d_1^{-1}(\bo)$ and $d_2^{-1}(\bo)$, so that the local
topological degrees $\deg_{\bo}(d_0)$, $\deg_{\bo}(d_1)$ and $\deg_{\bo}(d_2)$ are defined.

Then both $J(0,x)$ and $J(t,x)$ have an isolated critical point at the origin. If $0\neq t$ is
sufficiently close to zero then
$$\chi(\tilde{M}_t^{-})=\chi(M_t^{-})=1-(\deg_{\bo}(d_0)+\deg_{\bo}(d_1)+\operatorname{sign}(t)\cdot \deg_{\bo}(d_2))/2 ,$$
$$\chi(\tilde{L}_0)=\chi(L_0)=2\cdot (  1-\deg_{\bo}(d_0  )     ).$$
In particular, if $n=2$ then $\tilde{L}_0$ is finite and $\#\tilde{L}_0=2\cdot (  1-\deg_{\bo}(d_0  )     )$.
\end{theorem}

It is proper to add that there exists an efficient computer program which may compute the local topological degree
(see \cite{leckiszafraniec1}).

\section{Number of half--branches}\label{liczbagalezi}

In this section we shall show how to adopt some techniques developed
in \cite{nowelszafraniec, szafraniec7,szafraniec14} so as to compute the number of half--branches
of an analytic set of dimension $\leq 1$ emanating from a singular point.

Let $\Opis_n=\R\{t,x_1,\ldots,x_{n-1}\}$ denote the ring of germs
at the origin of real analytic functions.
If $I$ is an ideal in $\Opis_n$, let $V(I)\subset\R^n$ denote
the germ of zeros of $I$ near the origin, and let
$V_{\C}(I)\subset\C^n$ denote the germ of complex zeros of $I$.

\begin{rem} \label{powtorkadim}
If $I$ is proper then
$\dim_{\R}\Opis_n/I<\infty$ if and only if
$V_{\C}(I)=\{\bO\}$.
\end{rem}

Let $w_1,\ldots,w_m\in\Opis_n$, where $m\geq n-1$, be germs vanishing at the origin.
We shall write $\langle w_1,\ldots,w_m\rangle$ for the ideal in $\Opis_n$ generated by $w_1,\ldots,w_m$.

Let $W\subset \Opis_n$ denote the ideal generated by $w_1,\ldots, w_m$
and all $(n-1)\times (n-1)$--minors of the Jacobian matrix
$[ \partial w_i/\partial x_j]$. The ideal $W$ is proper
if and only if the rank of this matrix at the origin is $\leq n-2$.

If $V(W)=\{\bo\}$ then by the implicite function theorem the germ $V(w_1,\ldots,w_m)$ is of dimension $\leq 1$,
so that this set is locally an union of a finite family of half-branches emanating from the origin.
We shall say that $V(w_1,\ldots, w_m)$ is a curve having an algebraically isolated singularity
at the origin if  $W$ is proper and $\dim_{\R} \Opis_n/W<\infty$. 

From now on we shall assume that $m=n=3$
Let $M(3,3)$ denote the space of all $3 \times 3$--matrices with coefficients in $\R$.
By  \cite[Theorem 3.8]{nowelszafraniec} and comments in \cite[p. 1012]{nowelszafraniec}  we have

\begin{theorem}\label{powtorka2}
Assume that $V(w_1,w_2,w_3)$ is a curve having an  algebraically isolated 
singularity. 
 There exists a proper algebraic subset
$\Sigma\subset M(3,3)$
 such that for every non-singular matrix
$[a_{sj}] \in M(3,3)\setminus\Sigma$ and
$g_s=a_{s,1}w_1+a_{s,2}w_2+a_{s,3}w_3$, where $1\leq s\leq 3$,
the  set $V(g_1,g_2)$ is a curve having an algebraically isolated singularity at the origin and
$V(w_1,w_2,w_3)=V(g_1,g_2,g_3)\subset V(g_1,g_2)$.

In particular, if $V(w_1,w_2)$ is a curve having an algebraically isolated singularity then
one may take $g_s=w_s$.

If that is the case and  $J_p =\langle g_1,g_2,g_{3}^p\rangle$, where $p=1,2$,
then $J_2\subset J_1$ and $\dim_{\R}(J_1/J_2)<\infty$.
\end{theorem}

From now on we shall assume that
\begin{equation}\label{powtorkaw1}
\dim_{\R}\Opis_3/\langle t,g_1,g_2\rangle <\infty.
\end{equation}
As $\dim_{\R}(J_1/J_2)<\infty$ and $g_3(\bo)=0$, then by the Nakayama lemma
$\xi=\min\{s\ |\ t^s\cdot g_3\in J_2\}$ is finite. 
(In \cite{szafraniec14} there are presented  effective methods for computing this number.)
Let $k>\xi$ be an even positive integer.

Now we shall adopt to our case some arguments presented in \cite[pp. 529-531]{szafraniec14}.
There are germs $h_1,h_2,h_3\in\Opis_3$ such that
$$t^{\xi}g_3=h_1g_1+h_2g_2+h_3 g_3^2.$$
Let $Y_{\C}=V_{\C}(g_1,g_2)\setminus V_{\C}(g_3)$.
By (\ref{powtorkaw1}), the germ $t^k$ does not vanish at points in $V_{\C}(g_1,g_2)\setminus\{\bo\}$.
If $(t,x_1,x_2)=(t,x)\in Y_{\C}$ lies sufficiently close to the origin then $|h_3(t,x)|<M$ for some $M>0$,
$g_1(t,x)=g_2(t,x)=0$ and $g_3(t,x)\neq 0$. Hence
$$|g_3(t,x)|\geq |t|^{\xi}/M>|t|^k.$$
Then the origin is isolated in both $V_{\C}(g_3\pm t^{k},g_1,g_2)$.

Take $(t,x)\in V(g_1,g_2)\setminus\{\bo\}$ near the origin. By (\ref{powtorkaw1}), $t\neq 0$.
If $g_3(t,x)\neq 0$ then $g_3(t,x)\pm t^k$ has the same sign as $g_3(t,x)$.
If $g_3(t,x)=0$ then $g_3(t,x)+t^k>0$ and $g_3(t,x)-t^k<0$.
Write $b_+$ (resp. $b_-$, $b_0$) for the number of half-branches of $V(g_1,g_2)$
on which $g_3$ is positive (resp. $g_3$ is negative, $g_3$ vanishes). Put
$$H_{\pm}=\left(\frac{\partial(g_3\pm t^k,g_1,g_2)}{\partial(t,x_1,x_2)},g_1,g_2\right):
\R^3,\bo\rightarrow\R^3,\bo .$$
By \cite[Theorem 3.1]{szafraniec7} or \cite[Theorem 2.3]{szafraniec14}, the origin is isolated
in both $H_{\pm}^{-1}(\bo)$ and 
$$b_{+}+b_0-b_{-}=2\, \deg_{\bo}(H_+),$$
$$b_{+}-b_0-b_{-} =2\, \deg_{\bo}(H_-).$$

\begin{theorem}\label{powtorka3}
If $\dim_{\R}\Opis_3/\langle t,g_1,g_2\rangle<\infty$ then
the number $b_0$ of half-branches of $V(w_1,w_2,w_3)$ emanating
from the origin equals $\deg_{\bo}(H_+)-\deg_{\bo}(H_-)$.
\end{theorem}
\noindent{\em Proof.} As the matrix $[a_{sj}]$ is non-singular, then
$V(w_1,w_2,w_3)=V(g_1,g_2,g_3)$. Of course, $b_0$ equals the number
of half-branches of $V(g_1,g_2,g_3)$. Moreover,
$$b_0=\frac{1}{2}((b_++b_0-b_1)-(b_+-b_0-b_1))=\deg_{\bo}(H_+)-\deg_{\bo}(H_-).\ \Box$$

Now we shall explain how to compute the number of half-branches of $V(w_1,w_2,w_3)$
in the region where $t>0$.

\begin{prop}\label{powtorka4}
Put $g_i'(t,x)=g_i(t^2,x)$. Then $\dim_{\R}\Opis_3/\langle t,g_1',g_2'\rangle<\infty$
and $V(g_1',g_2')$ has an isolated singularity at the origin.
\end{prop}
\noindent{\em Proof.} By (\ref{powtorkaw1}), as $V_{\C}(t,g_1,g_2)=\{\bo\}$
then $V_{\C}(t,g_1',g_2')=\{\bo\}$. By Remark \ref{powtorkadim},
$\dim_{\R}\Opis_3/\langle t,g_1',g_2'\rangle<\infty$. We have
$$\frac{\partial(g_i',g_j')}{\partial(t,x_p)}(t,x)=2t\frac{\partial(g_i,g_j)}{\partial(t,x_p)}(t^2,x),\ \ \ 
\frac{\partial(g_i',g_j')}{\partial(x_1,x_2)}(t,x)=\frac{\partial(g_i,g_j)}{\partial(x_1,x_2)}(t^2,x),$$
and then $V(g_1',g_2')$ is a curve having an algebraically isolated singularity at the origin. $\Box$

\begin{rem}\label{powtorka5}
Let $J_p'=\langle g_1',g_2',(g_3')^p\rangle$. 
Put $\xi'=\min\{s\ |\ t^s\cdot g_3'\subset J_2'\}$. Of course, $\xi'\leq 2\cdot\xi$.
\end{rem}

Applying the same methods as above, one may compute the number $b_0'$ of half-branches
of $V(g_1',g_2',g_3')$. Obviously $b_0'/2$ equals the number of half-branches
of $V(w_1,w_2,w_3)$ lying in the region where $t>0$.

Other methods of computing the number of half-branches were presented in \cite{aokietal1}, \cite{aokietal2}
\cite{cuckeretal}, \cite{damon1}, \cite{damon2}, \cite{fukudaetal}, \cite{montaldivanstraten}.

According to  Khimshiashvili \cite{khimshiashvili1, khimshiashvili2},
if a germ $f:\R^2,\bo\rightarrow \R,0$ has an isolated critical point at the origin then
the number of real half--branches in $f^{-1}(0)$
equals $2\cdot(1-\deg_{\bO}(\nabla f))$, where $\nabla f:\R^2,\bO\rightarrow \R^2,\bO$ is the gradient
of $f$.

\section{Mappings between curves}\label{owale}
In this section we give sufficient  conditions for a mapping between
some smooth plane curves to have only  non-degenerate critical points.

Let $f=(f_1,f_2):\R^2\rightarrow\R^2$ be a smooth mapping.
Put $g=f_1^2+f_2^2$. Assume that $\delta^2>0$ is a regular value of $g$ and
$P=g^{-1}(\delta^2)$ is non-empty, so that $P$ is a smooth curve.
Obviously, $P=f^{-1}(S^1(\delta))$ and
$f|P:P\rightarrow S^1(\delta)$ is a smooth mapping between 1-dimensional manifolds.

At any $p\in P$ the gradient $\nabla g(p)=(\partial g/\partial x_1(p),\partial g/\partial x_2(p))$
is a non-zero vector perpendicular to $P$, and the vector
$T(p)=(-\partial g/\partial x_2(p),\partial g/\partial x_1(p))$
 obtained by rotating $\nabla g(p)$ counterclockwise
by an angle of $\pi/2$ is tangent to $P$. This way $T:P\rightarrow\R^2$ is a non-vanishing tangent
vector field along $P$.

Take $p\in P$. There exists a smooth maping 
$x(t)=(x_1(t),x_2(t)):\R\rightarrow P$ such that $x(0)=p$ and 
$x'(t)=T(x(t))$.  Hence
\begin{equation} \label{rowale1}
x_1'(t) = \left. -2\cdot\left( f_1\frac{\partial f_1}{\partial x_2} +   f_2\frac{\partial f_2}{\partial x_2}\right)\right|_{ (x(t))},\end{equation}
$$
 x_2'(t)= \left. 2\cdot\left( f_1\frac{\partial f_1}{\partial x_1} +  f_2\frac{\partial f_2}{\partial x_1}\right)\right|_{ (x(t))}  .
$$

As $g(x(t))=\delta^2$, then
$f(x(t))=(\delta\cos\theta(t),\delta\sin\theta(t))$ for some smooth function
$\theta:\R,0\rightarrow\R$.
Of course, $(\delta\cos\theta(0),\delta\sin\theta(0))=f(x(0))=f(p)$.
Applying the complex numbers notation we may write
\begin{equation} \delta\cdot e^{i\theta}=f_1(x(t))+i\, f_2(x(t)),\mbox{ where }  i=\sqrt{-1}. \label{rowale2}
\end{equation}

Put $J=\partial(f_1,f_2)/\partial(x_1,x_2)$ and $F_j=\partial(f_j,J)/\partial(x_1,x_2)$, where $j=1,2$.

\begin{lemma}\label{owale1}
A point $p\in P$ is a critical point of $f|P:P\rightarrow S^1(\delta)$ if and only if $J(p)=0$.

\end{lemma}
\noindent {\em Proof.} By (\ref{rowale1}), the derivative of the equation (\ref{rowale2}) equals
$$i\, \delta\, \theta'\cdot e^{i\theta}=
\left( \frac{\partial f_1}{\partial x_1}\, x_1'+\frac{\partial f_1}{\partial x_2}\, x_2'       \right)
+i\cdot \left( \frac{\partial f_2}{\partial x_1}\, x_1'+\frac{\partial f_2}{\partial x_2}\, x_2'       \right) $$
$$=2i(f_1+i\, f_2)\cdot J=2 i \,\delta\cdot e^{i\theta}\cdot J.$$
So $p\in P$ is a critical point of $f|P$ if and only if  ${\theta}'(0)=0$, i.e. if $J(p)=0$.  $\Box$

\begin{lemma}\label{owale2}
Suppose that $p\in P$ is a critical point of $f|P:P\rightarrow S^1(\delta)$. Then
$$\operatorname{sign}\,(\theta''(0))=\left. \operatorname{sign}\,\left( f_1\cdot F_1  +
f_2\cdot F_2     \right)\right|_{p}.$$
In particular, a point $p\in P$ is a non-degenerate critical point of $f|P:P\rightarrow S^1(\delta)$
if and only if $J(p)=0$ and $\left. ( f_1\cdot F_1  +
f_2\cdot F_2   )\right|_{p}\neq 0$.
\end{lemma}
\noindent {\em Proof.} Since $\theta'(0)=0$ and $J(p)=0$, after computing the second derivative of (\ref{rowale2})
the same way as above one gets
$$\left. i\,\delta\,\theta''\cdot e^{i\theta}\right|_0=\left. 2i\,\delta\cdot e^{i\theta}\cdot\left( \frac{\partial J}{\partial x_1}\, x_1'+\frac{\partial J}{\partial x_2}\,x_2'          \right)\right|_0$$
$$=\, 4\, i\,\delta\cdot e^{i \theta(0)}\cdot\left. \left(f_1\cdot F_1+
f_2\cdot F_2   \right)\right|_p.\ \Box$$

\begin{lemma}\label{owale3}
Let $f=(f_1,f_1):\R^2,\bo\rightarrow\R^2,\bo$ be an analytic mapping such that $J(\bo)=0$, and the origin is isolated
in both $f^{-1}(\bo)$ and $\nabla J^{-1}(\bo)$.

If $0<\delta\ll r\ll 1$ then $\tilde{S}^1(\delta)=D(r)\cap f^{-1}(S^1(\delta))$ is diffeomorphic to a circle,
$\tilde{D}^2(\delta)=D(r)\cap f^{-1}(D^2(\delta))$ is diffeomorphic to a disc, and
$f:\tilde{S}^1(\delta)\rightarrow S^1(\delta)$ has only non-degenerate critical points.
Moreover the one-dimensional set $J^{-1}(0)$ consisting of of critical points of $f$ is transverse to $\tilde{S}^1(\delta)$.
\end{lemma}
\noindent {\em Proof.} If the origin is isolated in $J^{-1}(0)$ then $f|\R^2\setminus\{\bo\}$ is a submersion
near the origin, and so $f:\tilde{S}^1(\delta)\rightarrow S^1(\delta)$ has no critical
points.

In the other case, $J^{-1}(0)\setminus\{\bo\}$ is locally a finite union of analytic half-branches emanating from the origin.
Let $B$ be one of them. The gradient $\nabla J(p)$ is a non-zero vector perpendicular to $T_p B$ at any $p\in B$.

The origin is isolated in $f^{-1}(\bo)$. 
By the curve selection lemma one may assume that
 $(f_1^2+f_2^2)|B$ has no critical points, so that $\nabla J$ and 
$$\nabla(f_1^2+f_2^2)=\left(2 f_1\frac{\partial f_1}{\partial x_1} + 2 f_2\frac{\partial f_2}{\partial x_1},
2 f_1\frac{\partial f_1}{\partial x_2} + 2 f_2\frac{\partial f_2}{\partial x_2}            \right)$$
are linearly independent along $B$. Then
$$0\neq \nabla J\times\nabla(f_1^2+f_2^2)=2 f_1\frac{\partial(J,f_1)}{\partial (x_1,x_2)} +
2 f_2\frac{\partial(J,f_2)}{\partial (x_1,x_2)}=-2(f_1\cdot F_1+f_2\cdot F_2)$$
along $B$. By previous lemmas, $f:\tilde{S}^1(\delta)\rightarrow S^1(\delta)$ has only
non-degenerate critical points. Other assertions are rather obvious. $\Box$

\section{Families of self-maps of $\R^2$}\label{rodziny}
In this section we  investigate 1-parameter families of plane-to-plane analytic mappings 

Let $f=(f_1,f_2):\R\times \R^2,\bo\rightarrow\R^2,\bo$ be an analytic function
defined in a neighbourhood of the origin. We shall write $f_t(x_1,x_2)=f(t,x_1,x_2)$
for $t$ near zero. Define three germs  $\R\times\R^2,\bo\rightarrow\R$ by
$$J=\frac{\partial(f_1,f_2)}{\partial(x_1,x_2)},\ 
F_i=\frac{\partial (f_i,J)}{\partial (x_1,x_2)}.$$
Put $J_t(x_1,x_2)=J(t,x_1,x_2)$.

From now on we shall also assume that 
\begin{equation}\label{rodzinyw1}
\begin{array}{l}
\dim_{\R}\Opis_3/\langle t,f_1,f_2\rangle<\infty, \ \dim_{\R}\Opis_3/\langle  t,F_1,F_2   \rangle <\infty,\\
J(\bo)=0,\ \dim_{\R}\Opis_3/\langle t,\frac{\partial J}{\partial x_1} , \frac{\partial J}{\partial x_2}      \rangle<\infty,
\end{array}
 \end{equation}
i.e. the origin is isolated in both $(\{0\}\times \C^2)\cap V_{\C}(f_1,f_2)$, $(\{0\}\times \C^2)\cap V_{\C}(F_1,F_2)$,
and $J_0$ has an algebraically isolated critical point at the origin.

\begin{lemma}\label{rodziny0}
Let $Q=\Opis_3/\langle t,J,F_1,F_2\rangle$. Then $\dim_{\R}Q<\infty$, i.e. the origin is isolated in
$(\{0\}\times\C^2)\cap V_{\C}(J,F_1,F_2)$.
\end{lemma}
\noindent{\em Proof.} Of course $\langle t,F_1,F_2\rangle\subset\langle t,J,F_1,F_2\rangle$.
Then $\dim_{\R}Q\leq\dim_{\R}\Opis_3/\langle t, F_1,F_2\rangle<\infty$. $\Box$\\[1em]

We shall write $g=f_1^2+f_2^2$ and $g_t(x_1,x_2)=g(t,x_1,x_2)$.
There exists a small $r_0>0$ such that $f_0^{-1}(\bo)\cap D^2(r_0)=\{\bo\}$.
For $|t|\ll \delta\ll r_0$, put $\tilde{S}_t^1(\delta)=f_t^{-1}(S^1(\delta))\cap D^2(r_0)$
and $\tilde{D}_t^2(\delta)=f_t^{-1}(D^2(\delta))\cap D^2(r_0)$.
If $\delta^2$ is a regular value of $g_0|D^2(r_0)$, then it is also a regular value of $g_t|D^2(r_0)$.
If that is the case then $\tilde{S}_t^1(\delta)$
is diffeomorphic to $\tilde{S}_0^1(\delta)\simeq S^1(1)$.
By the same argument,  $\tilde{D}_t^2(\delta)$
is diffeomorphic to $\tilde{D}_0^2(\delta)\simeq D^2(1)$.

By Lemmas \ref{owale2}, \ref{owale3} we get
\begin{lemma}\label{rodziny1}
 Critical points of $f_0:\tilde{S}_0^1(\delta)\rightarrow S^1(\delta)$
are non-degenerate, and 
$C(f_0|\tilde{S}_0^1(\delta))=\tilde{S}_0^1(\delta)\cap\{J_0=0\}$.

For $t$ near zero, critical points of  of $f_t:\tilde{S}_t^1(\delta)\rightarrow S^1(\delta)$
are non-degenerate too, and the number of critical points
$\#C(f_t|\tilde{S}_t^1(\delta))$  equals $\#(\tilde{S}_0^1(\delta)\cap \{J_0=0\})$.
Moreover the set of critical points of $f_t$, i.e. $J_t^{-1}(0)$, is transverse to $\tilde{S}_t^1(\delta)$. $\Box$
\end{lemma}

Let $I$ denote the ideal in the ring $\Opis_3$ 
generated by $J, F_1, F_2$, and let $V(I)\subset\R\times\R^2$ denote a representative
of the germ of zeros of $I$ near the origin. By Lemma \ref{rodziny0}, there exists $0<\delta\ll  1$
such that $\{0\}\times\tilde{D}_0^2(\delta)\cap V(I)=\{\bo\}$, and
$\{t\}\times \tilde{S}_t^1(\delta)\cap V(I)=\emptyset$ for $t$ sufficiently close to zero.
Put $\Sigma_t=\{x\in\tilde{D}_t^2(\delta)\, |\, (t,x)\in V(I)\}$. Hence $\Sigma_0=\{\bo\}$ and $\Sigma_t$
is contained in the interior of $\tilde{D}_t^2(\delta)$.

Let $I'$ denote the ideal in $\Opis_3$ generated by germs
$J$, $F_1$, $F_2$, $\partial(F_1,J)/\partial (x_1,x_2)$ and $\partial(F_2,J)/\partial (x_1,x_2)$.
Suppose that  $V(I')=\{\bo\}$.
Hence $\{t\}\times \tilde{D}^2(\delta)\cap V(I')$ is empty for $0\neq t$ close to zero.
By Proposition \ref{powierzchnie2} one gets

\begin{lemma}\label{rodziny2}
Suppose that $0<\delta\ll 1$ and $0\neq t$ is sufficiently close to zero.
Then the set of critical points of $f_t:\tilde{D}_t^2(\delta)\rightarrow D^2(\delta)$
consists of fold points, and a finite family $\Sigma_t$ of cusp points. $\Box$
\end{lemma}

\begin{rem}\label{rodziny2plus1}
By \cite[Theorem 3.1]{fukudaishikawa},
if $0\neq t$ is sufficiently close to zero then
$\#\Sigma_t\leq\dim_{\R}Q$ and $\#\Sigma_t=\dim_{\R}Q\bmod 2$.
\end{rem}

For $t\neq 0$ we shall write $\Sigma_t^{\pm}=\{x\in\Sigma_t\ |\ \mu_t(x)=\pm 1\}$,
where $\mu_t(x)$ is the local topological degree of $f_t$ at $x$. Put
$\operatorname{cusp\ deg}( f_t) =\sum_{x\in \Sigma_t}\, \mu_t(x) = \#\Sigma_t^+ -\#\Sigma_t^-$.
By Lemmas  \ref{owale3}, \ref{rodziny1}, \ref{rodziny2} and Theorem \ref{powierzchnie1} we get

\begin{prop}\label{rodziny3}
Suppose that $0<\delta\ll 1$, and $0\neq t$ is sufficiently close to zero.
Then \begin{itemize}
\item[(i)] the pair $(\tilde{D}_t^2(\delta),\tilde{S}_t^1(\delta))$ is diffeomorphic
to $(D^2(1),S^1(1))$, and $f_t:\tilde{D}_t^2(\delta)\rightarrow D^2(\delta)$
is such a mapping that $f_t^{-1}(S^1(\delta))=\tilde{S}_t^1(\delta)$,
\item[(ii)]  every point in $\tilde{D}_t^2(\delta)$ is either a fold point, a cusp point
or a regular point, and there is  a finite family of cusps
which all belong to $\tilde{D}_t^2(\delta)\setminus\tilde{S}_t^2(\delta)$,
\item[(iii)] $f_t|\tilde{S}_t^1:\tilde{S}_t^1(\delta)\rightarrow S_t^1(\delta)$ is locally stable,
and the set of critical points of $f_t$, i.e. $J_t^{-1}(0)$, is transverse to $\tilde{S}_t^1(\delta)$,
\item[(iv)] $\operatorname{cusp\  deg}(f_t)=2\chi(\tilde{M}_t^{-})+\deg(f_t|\tilde{S}_t^1(\delta))
-1-\#C(f_t|\tilde{S}_t^1(\delta))/2$\\
$$= 2\chi(\tilde{M}_t^{-})+\deg_0(f_0)-\#C(f_0|\tilde{S}_0^1(\delta)  ) /2-1,$$
\end{itemize}
where $\tilde{M}_t^{-}=\{ x\in \tilde{D}_t^2(\delta)\ |\ J_t(x)\leq 0 \}$. $\Box$
\end{prop}

Let $d_1,d_2:\R\times\R^2,\bo\rightarrow\R\times\R^2,\bo$
be defined as in Section \ref{kielki}. 

\begin{theorem}\label{rodziny4}
Let $f=(f_1,f_2):\R\times \R^2,\bo\rightarrow\R^2,\bo$ be an analytic function
defined in a neighbourhood of the origin such that $(\ref{rodzinyw1})$ holds.
 Suppose that the origin is isolated in $V(I')$, $d_1^{-1}(\bo)$ and $d_2^{-1}(\bo)$.

Then there exits $r>0$ such that the set of critical points of $f_t:D^2(r)\rightarrow \R^2$,
where $0\neq t$ is sufficiently close to zero,
consists of fold points, and a finite family $\Sigma_t$ of cusp points.
Moreover, the origin is isolated in $f_0^{-1}(\bo)$ and
$$\operatorname{cusp\ deg}(f_t)= \deg_{\bo}(f_0)-\deg_{\bo}(d_1)-\operatorname{sign}(t)\cdot\deg_{\bo}(d_2).$$

\end{theorem}
\noindent{\em Proof.} For any small $\delta>0$ there is $r>0$ such that
$D^2(r)\subset\tilde{D}_0^2(\delta)\setminus\tilde{S}_0^1(\delta)$,
so that also $D^2(r)\subset \tilde{D}_t^2(\delta)\setminus\tilde{S}_t^1(\delta)$
if $|t|$ is small.

By Lemma \ref{rodziny2}, the set of critical points of $f_t|\tilde{D}_t^2(\delta)$
consists of fold points, and a finite family $\Sigma_t$ of cusp points. Because $\Sigma_0=\{\bo\}$
then  $\Sigma_t$ is the set of cusp points of $f_t|D^2(r)$.

By (\ref{rodzinyw1}), the germ
$d_0=\nabla J_0:\R^2,\bo\rightarrow\R^2,\bo$ has an isolated zero at the origin.
By Theorem \ref{powtorka1} and Lemma \ref{rodziny1},
$$\#C(f_t|\tilde{S}_t^1(\delta))=\#(\tilde{S}_0^1(\delta)\cap\{J_0=0\})=2\cdot(1-\deg_{\bo}(d_0)),$$
for $0\neq t$ sufficiently close to zero. Our assertion is then a consequence of Proposition \ref{rodziny3}
and Theorem \ref{powtorka1}. $\Box$

Put $J'=J(t^2,x_1,x_2)$, $F_i'=F_i(t^2,x_1,x_2)$.

\begin{lemma}\label{rodziny5}
Suppose that $V(I')=\{\bo\}$. Then $\dim V(J,F_1,F_2)\leq 1$ and $\dim V(J',F_1',F_2')\leq 1$.

Moreover, if  $\dim_{\R}\Opis_3/I'<\infty$ then $V(J',F_1',F_2')$, as well as $V(J,F_1,F_2)$,
is a curve having an algebraically isolated singularity.
\end{lemma}
\noindent{\em Proof.} We have
$$\{\bo\}=V(I')=V(J,F_1,F_2)\cap V\left( \frac{\partial (F_1,J)}{\partial (x_1,x_2)} , \frac{\partial (F_2,J)}{\partial (x_1,x_2)}   \right),$$
so by the implicite function theorem $\dim V(J,F_1,F_2)\leq 1$. Of course,
$(t,x_1,x_2)\in V(J',F_1',F_2')$ if and only if $(t^2,x_1,x_2)\in V(J,F_1,F_2)$.
Hence $\dim V(J',F_1',F_2')\leq 1$ too.

 The ideal
$$K=\left\langle   J',F_1',F_2', \frac{\partial(F_1',J')}{\partial(x_1,x_2)} , \frac{\partial(F_2',J')}{\partial(x_1,x_2)}      \right\rangle\subset\Opis_3$$
is contained in the ideal $L$ generated by $J',F_1',F_2'$ and all $2\times 2$-minors of the derivative matrix of $(J',F_1',F_2')$.

As $\dim_{\R}\Opis_3/I'<\infty$, by the local Nullstellensatz, the origin is isolated in the set
of complex zeros of $I'$. Since
$$\frac{\partial(F_i',J')}{\partial(x_1,x_2)}(t,x_1,x_2)=\frac{\partial(F_1,J)}{\partial(x_1,x_2)}(t^2,x_1,x_2),$$
the origin is isolated in the set of complex zeros of $K$.
Hence $\dim_{\R}\Opis_3/L\leq\dim_{\R}\Opis_3/K<\infty$, and then $V(J',F_1',F_2')$ is a curve
having an algebraically isolated singularity at the origin.
The proof of the last assertion is similar. $\Box$\\[1em]

Suppose that the origin is isolated in $V(I')$. Let $b_0$  (resp. $b_0'$) be the number of half branches
in $V(J,F_1,F_2)$ (resp. $V(J',F_1',F_2')$) emanating from the origin.

By Lemma \ref{rodziny0}, no half-branch is contained in $\{0\}\times\R^2$.
Then by the curve selection lemma the family of half-branches is a finite union
of graphs of continuous functions $t\mapsto x^i(t)\in\R^2$, where $t$ belongs either to $(-\epsilon,0]$
or to $[0,\epsilon)$, $0<\epsilon\ll 1$, $x^i(0)=\bo$, $1\leq i\leq b_0$ (resp. $1\leq i\leq b_0'$),
and those graphs meet only at the origin.

Hence, if $0<t\ll 1$ then
$$b_0=\#\Sigma_t+\#\Sigma_{-t}=\#\Sigma_t^+ +\#\Sigma_t^-+\#\Sigma_{-t}^+ +\Sigma_{-t}^- ,$$
$$b_0'/2=\#\Sigma_t=\#\Sigma_t^+  + \#\Sigma_t^- .$$
By Theorem \ref{rodziny4}, we have
$$ \deg_{\bo}(f_0)-\deg_{\bo}(d_1)-\deg_{\bo}(d_2) = \#\Sigma_t^+  -  \#\Sigma_t^- ,$$
$$  \deg_{\bo}(f_0)-\deg_{\bo}(d_1)+\deg_{\bo}(d_2) = \#\Sigma_{-t}^+ -\#\Sigma_{-t}^- .$$
Then we have

\begin{theorem}\label{rodziny6}
Suppose that assumptions of Theorem \ref{rodziny4} hold.
Then numbers $\#\Sigma_{\pm t}^{\pm}$, where $t>0$ is small, are determined by $b_0,b_0',\deg_{\bo}(f_0),\deg_{\bo}(d_1),\deg_{\bo}(d_2)$.

Moreover, if $\dim\Opis_3/I'<\infty$ then $V(J,F_1,F_2)$ and $V(J',F_1',F_2')$
are curves having an algebraically isolated singularity at the origin. In that case
one may apply Theorem \ref{powtorka3} so as to compute $b_0$ and $b_0'$.
In particular, if $dim_{\R}\Opis_3/I''<\infty$, where
$$I''=\left\langle   F_1,F_2,\frac{\partial(F_1,F_1)}{\partial(t,x_1)} ,
\frac{\partial(F_1,F_1)}{\partial(t,x_2)}, \frac{\partial(F_1,F_1)}{\partial(x_1,x_2)}   \right\rangle,$$
then $V(F_1,F_2)$ is a curve having an algebraically isolated singularity at the origin.
In that case one may take $g_1=F_1$, $g_2=F_2$, $g_3=J$.
\end{theorem}

\section{Examples}\label{przyklad}

Examples presented in this section were calculated with the help 
of {\sc Singular} \cite{singular} and  the computer program 
written by Andrzej {\L}{\c e}cki \cite{leckiszafraniec1}.\\[1em]
{\bf Example 1.} Let $f=(f_1,f_2)=(x_1^3+x_2^2+tx_1,x_1 x_2)$.
Since $\dim_{\R}\Opis_3/\langle t, f_1,f_2\rangle=5,$
$\dim_{\R}\Opis_3/\langle t, F_1,F_2 \rangle=7,$
$\dim_{\R}\Opis_3/\langle t,\frac{\partial J}{\partial x_1}, \frac{ \partial J}{\partial x_2}\rangle=2,$
then (\ref{rodzinyw1}) holds. Moreover,
$\dim_{\R}\Opis_3/I'=8$, $\dim_{\R}\Opis_3/\langle \frac{\partial J}{\partial t},\frac{\partial J}{\partial x_1},
\frac{\partial J}{\partial x_2}\rangle=1$, and
$\dim_{\R}\Opis_3/\langle J,\frac{\partial J}{\partial x_1},\frac{\partial J}{\partial x_2}\rangle=3$.
Then the origin is isolated in $V(I')$, $d_1^{-1}(\bo)$ and $d_2^{-1}(\bo)$.
Using the computer program by {\L}{\c e}cki one may compute
$\deg_{\bo}(f_0)=-1$, $\deg_{\bo}(d_1)=+1$ and $\deg_{\bo}(d_2)=-1$.
By Theorem \ref{rodziny4}, $\operatorname{cusp\, deg}(f_t)=\operatorname{sign}(t)-2$ for $0\neq t$ sufficiently close to zero.

By Lemma \ref{rodziny5}, the set $V(J,F_1,F_2)$, as well as $V(J',F_1',F_2')$, is a curve having an algebraically
isolated singularity. Hence we may appy techniques presented in Section \ref{liczbagalezi} so as
to compute the number of half-branches of those curves.

One may verify that
$\dim_{\R}\Opis_3/I''= 8$,
so that $V(F_1,F_2)$ is a curve with an algebraically isolated singularity at the origin.

Put $J_p=\langle F_1,F_2,J^p\rangle$, where $p=1,2$. 
In that case  $\xi=2$, and so $k=4$. As $\dim_{\R}\Opis_3/\langle t,F_1,F_2\rangle<\infty$,
then (\ref{powtorkaw1}) holds. Set
$$H_{\pm}=\left(\frac{\partial(J\pm t^4,F_1,F_2 )}{\partial(t,x_1,x_2)},F_1,F_2\right):\R^3,\bo\rightarrow\R^3,\bo.$$

One may compute $\deg_{\bo}(H_+)=+2$, $\deg_{\bo}(H_-)=-2$.
By Theorem \ref{powtorka3}, $V(J,F_1,F_2)$ is an union of four half-branches
emanating from the origin, i.e. $b_0=4$.

Now we shall apply the same techniques so as to compute the number of half-branches
of $V(J',F_1',F_2')$.
By Proposition \ref{powtorka4}, $V(F_1',F_2')$ is a curve with an algebraically isolated singularity at the origin.
Put $J_p'=\langle F_1',F_2',(J')^p\rangle$, where $p=1,2$.  By Remark \ref{powtorka5},
$\xi'\leq 4$ and so one may take $k=6$. 
Let
$$H_{\pm}'=\left(\frac{\partial(J'\pm t^6,F_1',F_2' )}{\partial(t,x_1,x_2)},F_1',F_2'\right):\R^3,\bo\rightarrow\R^3,\bo.$$
One may compute $\deg_{\bo}(H_+')=+1$, $\deg_{\bo}(H_-')=-1$.
Then $V(J',F_1',F_2')$ is an union of two half-branches
emanating from the origin, i.e. $b_0'/2=1$. Hence, if $0<t\ll 1$ then
$\#\Sigma_t^+=0$, $\#\Sigma_t^-=1$, $\#\Sigma_{-t}^+=0$ and $\#\Sigma_{-t}^-=3$.\\[1em]
{\bf Example 2.} Let $f=(f_1,f_2)=(x_1^4+x_2^4+x_1^2 x_2^2+tx_1,x_1 x_2+tx_2)$.
In that case
$\dim_{\R}\Opis_3/\langle t, f_1,f_2\rangle=8,$
$\dim_{\R}\Opis_3/\langle t, F_1,F_2 \rangle=24,$
$\dim_{\R}\Opis_3/\langle t,\frac{\partial J}{\partial x_1}, \frac{ \partial J}{\partial x_2}\rangle=9,$
$\dim_{\R}\Opis_3/I'=33$, $\dim_{\R}\Opis_3/\langle \frac{\partial J}{\partial t},\frac{\partial J}{\partial x_1},
\frac{\partial J}{\partial x_2}\rangle=3$, and
$\dim_{\R}\Opis_3/\langle J,\frac{\partial J}{\partial x_1},\frac{\partial J}{\partial x_2}\rangle=12$.
Then the origin is isolated in $V(I')$, $d_1^{-1}(\bo)$ and $d_2^{-1}(\bo)$.
One may compute  $\deg_{\bo}(f_0)=0$, $\deg_{\bo}(d_1)=+1$ and $\deg_{\bo}(d_2)=0$.
By Theorem \ref{rodziny4}, $\operatorname{cusp\, deg}(f_t)=-1$ for $0\neq t$ sufficiently close to zero,
i.e. $\#\Sigma_t^+-\#\Sigma_t^-=-1$.

As $\dim_{\R}\Opis_3/I''=45$ then $V(F_1,F_2)$ is a curve having an isolated singularity at the origin.
Let $J_p$ be defined the same way as in the previous example. One may verify that
$\xi=2$, and so $k=4$.  Put
$$H_{\pm}=\left(\frac{\partial(J\pm t^4,F_1,F_2 )}{\partial(t,x_1,x_2)},F_1,F_2\right):\R^3,\bo\rightarrow\R^3,\bo.$$
One may compute $\deg_{\bo}(H_+)=0$, $\deg_{\bo}(H_-)=-2$.
Then $V(J,F_1,F_2)$ is an union of two half-branches
emanating from the origin, i.e. $b_0=2$.

Because $f_t(x_1,x_2)=f_{-t}(-x_1,-x_2)$, then $b_0'/2=1$ and $\#\Sigma_t^+=\#\Sigma_{-t}^+$, $\#\Sigma_t^-=\#\Sigma_{-t}^-$. So in this case there is no need to compute $\deg_{\bo}(H_{\pm}')$.
Hence, if $t>0$ then
$\#\Sigma_t^+=\#\Sigma_{-t}^+=0$ and
$\#\Sigma_{t}^-=\#\Sigma_{-t}^-=1$.


Zbigniew SZAFRANIEC\\
Institute of Mathematics, University of Gda\'nsk\\
80-952 Gda\'nsk, Wita Stwosza 57, Poland\\
Zbigniew.Szafraniec@mat.ug.edu.pl\\

\begin{thebibliography}{99}


  
\bibitem{aokietal1} K. Aoki, T. Fukuda, T. Nishimura, On the number
   of branches of the zero locus of a map germ $(\R^n,0)\rightarrow
   (\R^{n-1},0)$. In {\em Topology and Computer Science: Proceedings of
   the Symposium held in honour of S.~Kinoshita, H.~Noguchi and
   T.~Homma on the occasion of their sixtieth birthdays\/},
   1987, pp. 347-363.
\bibitem{aokietal2} K. Aoki, T. Fukuda, T. Nishimura, An algebraic
   formula for the topological types of one parameter bifurcations
   diagrams. {\em Archive for Rational Mechanics and Analysis\/}
   {\bf 108} (1989), 247-265.






\bibitem{cuckeretal} F. Cucker, L. M. Pardo, M. Raimondo, T. Recio, M.-F. Roy,
On the computation of the local and global analytic branches of a real algebraic curve.
In {\em Applied Algebra, Algebraic Algorithms and Error-Corecting Codes}, Lecture Notes in
Computer Sci. {\bf 356} (Springer-Verlag, 1989), pp. 161-181.
\bibitem{damon1} J. Damon, On the number of branches for real and complex
    weighted homogeneous curve singularities. {\em Topology\/} {\bf 30}
    (1991), 223-229.
\bibitem{damon2} J. Damon, $G$-signature, $G$-degree, and symmetries of the
   branches of curve singularities. {\em Topology\/} {\bf 30} (1991),
   565-590.
  \bibitem{eisenbudlevine} D. Eisenbud, H. I. Levine, An algebraic formula for
      the degree of a \( C^{\infty} \)-map germ.
      {\em Annals of Mathematics\/},
   {\bf 106} (1977), 19-44.
\bibitem{singular}
G.-M.~Greuel, G.~Pfister, and H.~Sch\"onemann,
\newblock {{\sc Singular} 3.0.2}. A Computer Algebra System for
Polynomial Computations.
\newblock Centre for Computer Algebra, University of Kaiserslautern (2006).
\newblock {\tt http://www.singular.uni-kl.de},

\bibitem{fukudaetal} T. Fukuda, K. Aoki, W.Z. Sun, On the number of branches
        of a plane curve germ. {\em Kodai Math. Journal\/} {\bf 9} (1986),
     179-187.
\bibitem{fukudaishikawa} T. Fukuda, G. Ishikawa, On the number of cusps of stable
perturbations of a plane-to-plane singularity. {\em  Tokyo J.  Math.} {\bf 10}
(1987), 375-384.
\bibitem{fukuda2} T. Fukuda, Topological triviality of plane-to-plane  singularities.
Geometry and its applications (Yokohama, 1991), 29-37, World Sci. Publ., River Edge, NJ, 1993.
\bibitem{fukui2} T. Fukui, An algebraic fomula for a topological invariant
   of bifurcation of  1-parameter family of function-germs, {\em
   Stratifications, singularities, and differential equations, II(Marseille,
   1990; Honolulu, HI, 1990)\/} Travaux en Cours {\bf 55}, Hermann, Paris
   1997, pp. 45-54.
\bibitem{gaffneymond} T. Gaffney, D. Mond, Cusps and double folds of germs of analytic mappings $\C^2\rightarrow\C^2$.
{\em J. London Math. Soc.} {\bf 43} (1991), 185-192.



     \bibitem{khimshiashvili1} G. M. Khimshiashvili, On  the local degree of
        a smooth mapping. {\em Comm. Acad. Sci. Georgian SSR\/} {\bf 85}(1977),
        309-311 (in Russian).
\bibitem{khimshiashvili2} G. M. Khimshiashvili, On the local degree of
        a smooth mapping. {\em Trudy  Tbilisi Math. Inst.} {\bf 64} (1980),
        105-124.
\bibitem{krzyzanowskaszafraniec}
I. Krzy\.zanowska, Z. Szafraniec, On polynomial mappings from the plane to the plane. {\em J.   Math. Soc. Japan} {\bf 66} (2014), 805-818.


 \bibitem{leckiszafraniec1} A. \L\c{e}cki, Z. Szafraniec, Applications 
	of the Eisenbud \& Levine's theorem to real algebraic geometry.
	{\em Computational Algebraic Geometry\/} Progr. in Math. {\bf 109},
	Birkh\"{a}user 1993, pp. 177-184.
\bibitem{levine} H. I. Levine, Mappings of manifolds into the plane. {\em Amer. J. Math.} {\bf 88} (1966), 357-365.








\bibitem{montaldivanstraten} J. Montaldi, D. van Straten, One-forms on
     singular curves and the topology of real curve singularities.
     {\em Topology\/} {\bf 29} (1990), 501-510.
\bibitem{moyaballesteros1} J. A. Moya--P\'erez, J. J. Nu\~{n}o-Ballesteros, The link of a finitely determined
map germ from $\R^2$ to $\R^2$. {\em J. Math. Soc. Japan} {\em 62} (2010), 1069-1092.
\bibitem{moyaballesteros2} J. A. Moya--P\'erez, J. J. Nu\~{n}o-Ballesteros, Topological triviality
of families of map germs from $\R^2$ to $\R^2$. {\em J. of Singularities} {\bf 6} (2012), 112-123.
\bibitem{nowelszafraniec} A. Nowel, Z. Szafraniec, On the number of branches of a real 
curve singularities. {\em   Bull. London Math. Soc.} {\bf 43} (2011), 1004-1020.

  
\bibitem{quine} J. R. Quine, A global theorem for singularities of maps between oriented 2-manifolds.
{\em Trans. Amer. Math. Soc.} {\bf 236} (1978), 307--314.
\bibitem{rieger} J. H. Rieger, Families of maps from the plane to the plane. {\em J. London Math. Soc.} {\bf 36} (1987), 351-369.
  
  \bibitem{szafraniec7} Z. Szafraniec, On the number of branches of a
     1-dimensional semianalytic set. {\em Kodai Math. Journal\/}
   {\bf 11} (1988), 78-85.
  \bibitem{szafraniec14} Z. Szafraniec, A formula for the number
      of branches of one-dimensional semianalytic sets. {\em Math.
     Proc. Cambridge Phil. Soc.\/} {\bf 112} (1992), 527-534.
\bibitem{thom} R. Thom, Les singularit\'es des applications diff\'erentiables. {\em Ann. Inst. Fourier, Grenoble} {\bf 6} (1955-1956), 43-87.

\bibitem{whitney} H.~Whitney, On singularities of mapping of Euclidean spaces. I. Mappings of the plane into the plane. {\em Annals of Mathematics} {\bf 62} (1955), 374-410.





\end{thebibliography}
\end{document}